

\baselineskip=14pt
\parskip=10pt
\def\halmos{\hbox{\vrule height0.15cm width0.01cm\vbox{\hrule height
  0.01cm width0.2cm \vskip0.15cm \hrule height 0.01cm width0.2cm}\vrule
  height0.15cm width 0.01cm}}

\magnification=\magstephalf

\def\1{{\overline{1}}}
\def\2{{\overline{2}}}
\parindent=0pt
\overfullrule=0in
\def\Tilde{\char126\relax}
\def\frac#1#2{{#1 \over #2}}
\centerline
{\bf A  One-Line Proof of Leversha's ``Quartet of Isogonal Conjugates" Theorem }
\bigskip
\centerline
{\it Shalosh B. EKHAD}

In [L], Gerry Leversha spent quite a few pages, using human ingenuity and lots of previous knowledge, to prove the following
elegant theorem (see [W] for the definition of {\it Isogonal Conjugate}).

{\bf Leversha's Quartet Theorem} ([L]): Let $ABCD$ be a quadrilateral which is not cyclic. Let $A^{*}$ be the isogonal conjugate of
 $A$ with respect to $\Delta BCD$, and similarly let $B^{*}$,  $C^{*}$,  $D^{*}$  the isogonal conjugates of
 $\Delta ADC$, $\Delta ABD$, $\Delta ABC$ respectively. Then $A$,$B$,$C$, and $D$ are the circumcenters of
 $\Delta B^*C^*D^*$, $\Delta A^*C^*D^*$, $\Delta A^*D^*D^*$, and $\Delta A^*C^*B^*$.

{\bf Proof}: Download the Maple package {\tt RENE.txt} freely available from Doron Zeilberger's website \hfill\break
{\tt https://sites.math.rutgers.edu/{}\Tilde zeilberg/PG/RENE.txt} \quad,\hfill\break
start {\tt Maple}, and type {\tt read `RENE.txt`:} followed by {\tt Leversha();}, 
and in {\bf one nano-second}, you would get that it is {\tt true}. The Maple code for proving Leversha's theorem is as follows:

{\tt
T1 := Te(m, n); T2 := Te(m1, n1); A := [0, 0]; B := T1[3]; C := T2[3]; \hfill\break
evalb(DeSq(IsogonalConjugate(m1, n1, B), A) = DeSq(IsogonalConjugate(m, n, C), A))
}

{\bf Explanation by Doron Zeilberger}

Without loss of generality (by translating, rotating, and shrinking) we can assume that $A=(0,0)$, $D=(1,0)$,
$$
tan(\frac{1}{2} \angle BAD)=m \quad, \quad
tan(\frac{1}{2} \angle ADB)=n \quad, \quad tan(\frac{1}{2} \angle CAD)=M
\quad, \quad tan(\frac{1}{2} \angle ADC)=N \quad,
$$
so now $\Delta ADB$ and $\Delta ADC$ are what are called (in {\tt RENE.txt}),  {\tt Te(m,n)} and  {\tt Te(M,N)}. Using
the straightforward {\tt IsogonalConjugate} {\it macro} we get $B^*$ and $C^*$ and we just check that the `distance-squared' ({\tt DeSq}), between
$A$ and $B^*$ and $A$ and $C^*$ are the same. Of course, this is also true about the distance between
$A$ and $B^*$ and $A$ and $D^*$, and by the {\it transitivity} of the $=$ relation, all three distances are equal and hence
indeed $A$ is the circumcenter of  $\Delta B^*C^*D^*$.

But using this `analytic' approach (rather than Leversha's {\it synthetic} approach) gives much more. 
{\tt LevershaRadius(m,n,M,N)} gives that the  {\bf exact} value of the circumradius of $\Delta B^*C^*D^*$ is
$$
{\frac { \left( N-n \right)  \left( Nn+1 \right)  \left( {m}^{2}+1 \right)  \left( {M}^{2}+1 \right) }
{ \left( ( MN-1 )  ( mn-1) +
 ( M+N ) ( n+m )  \right)  \left( MN ( n+m ) -mn ( M+N) +M+N-m-n \right) }} \quad .
$$

{\tt RENE.txt} can also tell you the exact locations of $B^*$ and $C^*$:
$$
B^* \, =\,
( \, {\frac { \left( Mm+M-m+1 \right)  \left( Mm-M+m+1 \right)  \left( Nn+1 \right)  \left( N-n \right) }
{ \left( MN(mn-1)+(M+N)(m+n)-mn+1 \right)  \left( MN(m+n)-mn(M+N)+M+N-m-n \right) }} \, ,
$$
$$
{\frac {2\, \left( Nn+1 \right)  \left( N-n \right)  \left( Mm+1 \right)  \left( M-m \right) }
{ \left( MN(mn-1)+(M+N)(m+n)-mn+1 \right)  \left( MN(m+n)-mn(M+N)+M+N-m-n \right) }} ) \quad,
$$
$$
C^* \, =\,
( \, {\frac { \left( Mm+M-m+1 \right)  \left( Mm-M+m+1 \right)  \left( Nn+1 \right)  \left( -n+N \right) }
{ \left( MN(mn-1)+(M+N)(m+n)-mn+1 \right)  \left( MN(m+n)-mn(M+N)+M+N-m-n \right) }} \, ,
$$
$$
-\,{\frac { 2\,\left( Nn+1 \right)  \left( -n+N \right)  \left( Mm+1 \right)  \left( -m+M \right) }
{ \left( MN(mn-1)+(M+N)(m+n)-mn+1 \right)  \left( MN(m+n)-mn(M+N)+M+N-m-n \right) }} ) \quad.
$$
So we have a nice {\it lagniappe}, we have proved something {\bf stronger} than the fact that  $B^*$ and $C^*$ are equidistant from $A$.
It turns out that $B^*$ and $C^*$ are {\bf mirror-images} of each other with respect to the line $AD$. 
Of course this is true for all the other five combinations.
\halmos 

\bigskip

{\bf References}

[L] Gerry Leversha, {\it A quartet of isogonal conjugates}, The Mathematical Gazette, {\bf 100}, issue 548, (July, 2016), 336-341,
(Note 100.23). Available from JSTOR.

[W] The Wikipedia Foundation, {\it Isogonal Conjugate} \hfill\break
{\tt https://en.wikipedia.org/wiki/Isogonal\_conjugate}

\bigskip
\hrule
\bigskip
Shalosh B. Ekhad, c/o D. Zeilberger, Department of Mathematics, Rutgers University (New Brunswick), Hill Center-Busch Campus, 110 Frelinghuysen
Rd., Piscataway, NJ 08854-8019, USA. \hfill\break
Email: {\tt ShaloshBEkhad at gmail dot com}   \quad .

\end